\documentclass[final,epsfig,amsfont]{article}
\usepackage{latexsym}
\usepackage[dvips]{graphicx}
\newcommand{\eeq}{\end{equation}}
\newcommand{\beq}{\begin{equation}}

\def\R{\ensuremath{\mathbf R}}
\def\N{\ensuremath{\mathbf N}}
\def\Q{\ensuremath{\mathbf Q}}

\def\p{\ensuremath{\mathbf P}}

\def\bna{\mathbf N}

\begin{document}
\title{Extreme value laws for fractal intensity functions in
dynamical systems: Minkowski analysis}
\author{
Giorgio Mantica and
Luca Perotti${}^{\dag}$, \\
Center for Non--linear and Complex Systems, \\
Dipartimento di Scienza ed Alta Tecnologia, \\ Universit\`a dell'Insubria,\\
via Valleggio 11, 22100
Como, Italy. \\ Also at I.N.F.N. sezione di Milano, CNISM unit\`a di Como and INDAM, GNFM. \\
${}^{\dag}$ Department of Physics, Texas Southern University, Houston, Texas 77004 USA.\\
}
\date{}
\maketitle
\begin{abstract}
Typically, in the dynamical theory of extremal events, the function that gauges the intensity of a phenomenon is assumed to be convex and maximal, or singular, at a single, or at most a finite collection of points in phase--space. In this paper we generalize this situation to fractal landscapes, {\em i.e.} intensity functions characterized by an uncountable set of singularities, located on a Cantor set. This reveals the dynamical r\^ole of classical quantities like the Minkowski dimension and content, whose definition we extend to account for singular continuous invariant measures. We also introduce the concept of extremely rare event, quantified by non--standard Minkowski constants and we study its consequences to extreme value statistics. Limit laws are derived from formal calculations and are verified by numerical experiments.
\end{abstract}




\begin{center}
{\em Dedicated to the memory of Joseph Ford, on the twentieth anniversary of his departure}
\end{center}

\section{Introduction}
\label{sec1}

The study of extreme events is a fundamental chapter in statistics and a pillar in its applications to natural phenomena. In recent years, this theory has been approached from a dynamical perspective \cite{col}. Rather than arising from a stochastic process, extreme events occur along the {\em deterministic} evolution of a system: see the review \cite{jorge1} and references therein. This is not surprising, for {chaotic} properties of dynamical systems render them equivalent to stochastic processes, in a precise technical sense \cite{alexeev,joe-phystoday}. This intuition guided the founders of modern non--linear dynamics like Joseph Ford, to whose memory we dedicate this paper, and is still extremely fruitful.

In this theory, one has been able to recover the three basic attraction laws of the statistical theory and to reproduce a wealth of its characteristics. This is highly desirable, since it is clear that natural phenomena follow  deterministic, albeit chaotic, evolution laws. In this endeavor, the ``conventional'' approach---let us call it so, for clarity---has been to focus on the events which occur when the system's evolution approaches a {\em single} point, or at most a finite collection of points, in a low dimensional phase space \cite{extremal}: the closer the approach, the more extreme the associated phenomenon. As a consequence, the function that gauges the strength of extreme events (we shall call this function the {\em intensity} in what follows) whether a distance, or a ``physical'' function \cite{physical,towards}, is rather simple.

In this paper we want to pursue a different approach: we try to construct abstract models that, albeit amenable of analysis, reproduce the expected complexity of a {\em physical landscape} of the intensity function. Consider in fact a realistic dynamical system modeling the evolution of a natural phenomenon. Typically, its motions  evolve in a many dimensional phase space, in which the intensity of a physical phenomenon is a non--convex function with many local maxima, almost degenerate, or even singularities. We may think of projecting this evolution on a manifold of reduced dimension via any of the well known techniques: surfaces of section, embedding, {\em et cetera}. Clearly, the graph of the intensity function so projected becomes a very complicated object, possibly displaying a fractal structure.

Alternatively, consider the minimization of a non convex, multi--variable function---a hard classical problem of computational mathematics. In \cite{complex} an algorithm has been proposed, by which a chaotic motion explores phase--space, being attracted by minima when they appear in its neighborhood.  This technique has been applied to inverse fractal problems \cite{complex} and to model tuning for oil well production data \cite{oil}.

Motions where close encounters with a scattered set of points have explosive significance may also evolve in real spaces: think of a robot ``sniffing'' the air in search of odor sources while cruising around \cite{snif}. These sources can be mines in a field. Probability of sensing a mine decays exponentially with distance; mine distribution cannot be assumed uniform and certainly/unfortunately more than one mine is to be expected. A further example is obtained by substituting the robot with a human and mines with infected individuals: spread of infections can be modeled this way \cite{theo}, so that extreme values of a distance function become directly related to the probability of contagion.

The above are just a few examples that motivate the investigation presented in this paper: we shall study dynamical extrema in low dimensional systems, when the intensity function has a complex, hierarchical structure. In so doing, while being closer to physics, we also introduce a new framework for mathematical investigation that reveals the {\em dynamical} r\^ole of objects of deep significance, like Minkowski dimension and content.

In summary, our work extends the conventional theory in some respects: first, as mentioned, rather than defining the intensity as the distance from a finite sets of points, we let it be a rather complicated function---that we also derive from the distance from uncountable Cantor sets. Second, we consider the interplay of such sets and functions with motions with singular continuous invariant measure. This leads us to define generalized Minkowski dimensions and contents. Finally, we find that, under certain circumstances, limit laws hold also for merely ergodic dynamical systems with no decay of correlations.

Previous related investigations with significant links to our results must be quoted, while underlying differences. In \cite{physical,towards,turche}, extreme value statistics for intensity functions taking their maximum value at a specific point of a fractal attractor have been studied. It has been shown that these statistics are linked to local dimensions of invariant measure on the attractor and can be used to reveal these quantities. Singular continuous invariant measures have been considered in \cite{nonsmooth,turche}. These works, as well as other previous investigations we are aware of, are concerned with {\em local} quantities. To the contrary, we pursue here a {\em global} approach. This difference emerges clearly when considering the relation between extreme value laws and hitting/return times statistics \cite{extremal}: conventional local quantities are replaced in our theory by global hitting/return times statistics, like those studied in \cite{prluno,miojstat}.
Finally, still in the local framework, extreme value laws have been shown to hold also in non--mixing systems, via the peak over threshold approach \cite{potpot}, while in this paper they are observed within the block maxima technique.

The approach adopted in this paper can be properly defined as experimental mathematics:
while the definitions of the objects under examination, as well as the form of the conjectured limit laws, are complete and rigorous and, frequently, the value of various constants involved are obtained from explicit formal computations, we do {\em not} provide proofs of the conjectured results. We only put forward a heuristic argument in Section \ref{sec-toy}, which might perhaps lead to rigorous analysis. In fact, to the best of our knowledge, the examples presented herein belong to classes not studied so far that might require the development of new techniques. We therefore bring our contribution by exposing this new framework and by providing numerical evidence.

Let us now outline the organization of this paper. In the next section we sketch the conventional steps of dynamical extreme value theory, mainly to fix notations, and we describe the set-up of our numerical experiments.

In the third section we introduce a simple geometrical model of intensity function with a hierarchical structure: the Cantor Ladder---not to be confused with the more common Cantor Staircase. 
Readers who want to skip this specific example may go directly to Section \ref{sec-discantor}.
We show that extreme value laws hold for this function, both for a chaotic dynamics and for a non-chaotic, just ergodic, dynamical system, with no decay of correlations. Even if the form of the extreme distributions can be found {\em exactly}, it cannot so far be proven {\em rigorously}, but only justified by numerical experiments. It therefore constitutes an interesting example for further investigation. The same remark also applies to the cases that follow.

In Section \ref{sec-discantor}, we show that the previous example can be related to a generalization of the conventional setting described above, in which the intensity of a phenomenon is a function of the distance from prescribed points in phase space. In fact, we let these points populate
a compact set $K$, the ternary Cantor set, and we show that the results obtained in the toy model of Section \ref{sec-toy} still hold.  A particular case is then put in evidence:
when the invariant measure of the dynamical system is the Lebesgue measure in $\R^n$ and the intensity function is the logarithm of the inverse distance to $K$,
extreme value statistics crucially depends on two quantities of fractal analysis: the Minkowski dimension and the Minkowski content of the set $K$. This is a remarkable situation in which these quantities take on a dynamical meaning.

In Section \ref{sec-genmink} we introduce a further element: in the same setting of Section \ref{sec-discantor}, we consider dynamical systems with a generic invariant measure $\mu$. In the conventional case, such measures have been studied in \cite{extremal,turche}, but still in relation to single--point distance functions.
In the proposed new framework, extreme value laws depend on {\em generalized} Minkowski dimension and content of the set $K$, {\em defined with respect to} the measure $\mu$. Although quite natural, we have not been able so far to locate in the literature previous mention of these generalized quantities. When they exist, the results of the previous section are reproduced in this generalization.

To the contrary, in Section \ref{sec-extrare}, we study a family of sets $K$, a dynamical map of M\"obius type and an invariant measure (also due to Minkowski, the remarkable Minkowski question--mark function) such that the conventional scaling law (which yields non--trivial dimension and content) does not hold.
A different law is derived and new Minkowski constants are introduced, that reflect the non--analytic behavior of the Minkowski question--mark function at rational values. Physically, we interpret this law as describing {\em extremely rare} events. Again, the new extreme value laws can be exactly computed and tested numerically with success.

A particularly interesting  instance of this non--analytic behavior is studied in Section \ref{sec-usuenne}: this is the case of a set $K$ whose Hausdorff and Minkowski dimension differ. We derive the form of the scaling law of extreme events and we test it numerically. Finally, the conclusion recap the highlights of this work in three paragraphs.

\section{A brief summary of extreme value theory}
\label{sec-extr}
In this section we briefly introduce the dynamical approach to extreme value theory, following the comprehensive review \cite{jorge1}, and we describe the set--up of our numerical experiments.

Consider a canonical dynamical system $(\chi,\varphi,\mu)$, comprising a transformation $\varphi:\chi \to \chi$ on a measure space $\chi$, which preserves a probability measure $\mu$ which is also ergodic. Also defined on $\chi$ is a physical observable
$f:\chi \to \R_+ \cup\{\infty\}$, that represents the {\em intensity} of the phenomenon under examination.
Dynamics $\varphi$ and observable $f$ permits to construct a stationary process $X_n$, defined as $X_n= f \circ \varphi^n$ for each $n\in \N$.
We let $F$ be the cumulative distribution function of the process:
\[F(s)=\mu \{ x \in \chi  \mbox{ s.t. } f(x) \leq s),
 \]
and we also let $\overline{F}=1-F$. Consider now the evolution starting from $x \in \chi$, over a finite time span of length $n$, and the maximum value reached by the intensity $f$ during this time interval. This defines the block maximum  $M_n(x)$:
\begin{equation}
\label{eq-mn}
M_n(x)=\max\{X_0(x),\ldots,X_{n-1}(x)\}.
\end{equation}
An \emph{Extreme Value Law} for $M_n$ exists if there is a non-degenerate distribution function  $G:\R\to[0,1]$, so that, for every $\tau>0$, there exists a sequence of levels $h_n=h_n(\tau)$, $n=1,2,\ldots$,  such that
\begin{equation}
\label{eq-un}
  n\mu \{ x \in \chi \mbox{ s.t. } f(x) >h_n \} \to \tau,\;\mbox{ as $n\to\infty$,}
\end{equation}
and for which the following holds:
\begin{equation}
\label{eq-evl}
 \mu \{ x \in \chi  \mbox{ s.t. } M_n(x) \leq h_n \} \to G(\tau),\;\mbox{ as $n\to\infty$}.
\end{equation}

The meaning of these definitions is appreciated in the case when $X_j$ are independent, equally distributed random variables with a probability measure $\p$. Clearly,
$
  \p(M_n < s)  = F(s)^n,
$
so that
$
  \log \p(M_n < s)  = n \log F(s) = n \log ( 1 - \bar{F} (s) ) \sim
  - n  \bar{F} (s)
$
when $s$ tends to the supremum of the support of the distribution of $f$ values. One then chooses a sequence of values $h_n$ such that
\begin{equation}
\label{eq-2}
  \tau_n := n \bar{F} (h_n) \to \tau, \;\mbox{ as $n\to\infty$}.
\end{equation}
Typically, $\tau$ is parameterized by a real variable $y$, $\tau = \tau(y)$, so that also $h_n$ depends on $y$. Then, in the limit of $n$ tending to infinity one has
\begin{equation}
\label{eq-3}
  \p(M_n < h_n(y)) = F(h_n(y))^n  \rightarrow e^{-\tau(y)}.
\end{equation}
When passing from independent random variables to dynamics, a similar law may still apply, with the introduction of an {\em extremal index} $\theta$ \cite{extremal}, so that
\begin{equation}
\label{eq-3bis}
  \mu( \{ x \mbox{   s.t.  } M_n(x) < h_n(y)\}) \rightarrow e^{-\theta \tau(y)}.
\end{equation}
The extremal index $\theta$ is brought about by dynamical correlations: in fact, at variance with the independent variables case, deterministic evolution may spend sizable time--spans in neighborhoods of points characterized by high values of the intensity function $f$, yielding a clustering of extrema. This phenomenon has been investigated in \cite{extremal} and in successive works by the same authors.

Finally, in the conventional theory it is required that an affine scaling of the function of the extrema $M_n$ be performed, with sequences of parameters $a_n$ and $b_n$, so that the extreme value law takes the form
\begin{equation}
\label{eq-ca502}
   \mu (\{ x \mbox{  s.t.  } a_n [{M}_n(x) - b_n ] \leq y \}) \rightarrow e^{-\tau_j(y)}
\end{equation}
where $\tau_j$ is one of three standard forms \cite{jorge1}.  To the contrary, in the following we will derive limiting laws that cannot be put in this form.

To verify the validity of eq. (\ref{eq-3bis}) and to infer dynamical parameters from observations two fundamental approaches have been proposed: block maxima and peaks over threshold. In the latter, the statistics of events whose intensity surpasses a certain value is investigated. This technique seems better suited for practical applications to experimental time series: see \cite{potpot} and references therein. To the contrary, in this paper we follow the former approach, that is, we compute experimentally the l.h.s. of eq. (\ref{eq-3bis}), by Birkhoff sums of numerically generated trajectories.  We compare it with the expected limit laws without estimating or fitting parameters, but using the theoretical values obtained otherwise.

The relation between the length of the Birkhoff sums ({\em i.e.} the total time of the evolution) and the length $n$ of the time windows upon which block maxima are computed has been investigated in \cite{faraturca}: we put ourselves in a situation so to have a sufficient number of samples $\{M_n(x_j), j=1,\ldots,J\}$, to estimate reliably the phase space average. When needed, to achieve this goal we perform numerical calculations on a cluster of processors. While this brute--force approach is simple to the border of trivial, a comment must be made on the results of Sections \ref{sec-genmink} -- \ref{sec-usuenne}, where a singular continuous, invariant measure is investigated.
In this case, direct evaluation of the dynamics $\varphi$ is numerically ill-conditioned and we employ the technique described in \cite{myprl}, \cite{prluno},  based on Iterated Function Systems.

\section{Fractal landscape: a Cantor Ladder}
\label{sec-toy}

\begin{figure}
\centerline{\includegraphics[width=.6\textwidth, angle = -90]{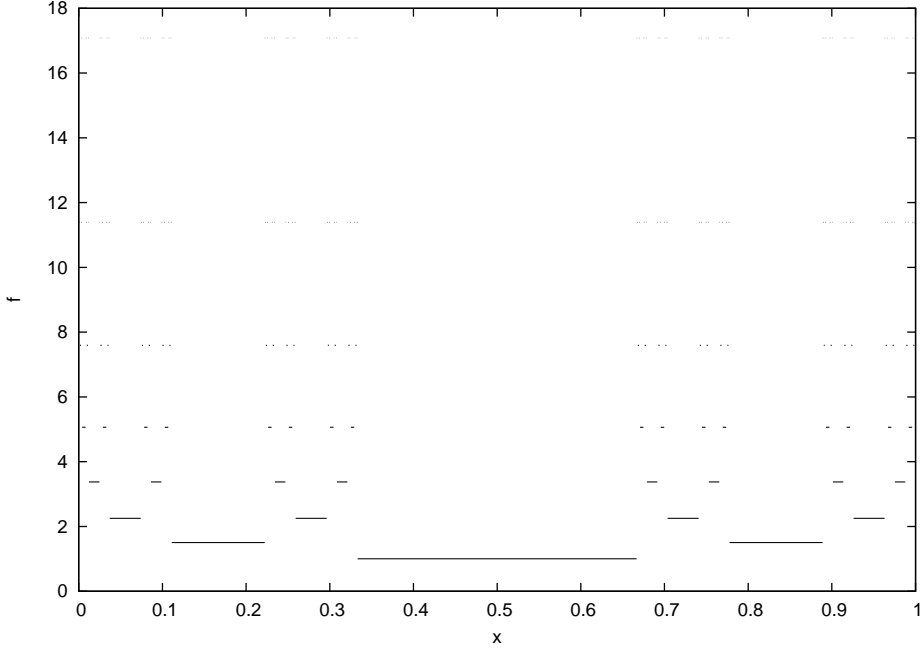}}
\caption{Intensity function $f$ taking the value $\beta^m$ on the gaps of order $m$ of the ternary Cantor set. Here, $\beta=3/2$. Only the lower part of the graph (that extends to infinity in the vertical direction) is shown.}
\label{disegno}
\end{figure}

We now define a dynamical model that exhibits extreme value features as observed in natural phenomena in a more suggestive way than what obtained by the ``conventional'' approach. Consider $\chi = [0,1]$ and the Cantor Ladder function depicted in Figure \ref{disegno}, which can be compared to the more common Cantor Staircase function.
This function takes constant values on every gap of the well known ternary Cantor set, of increasing magnitude as the size of the gap gets smaller. This generates a prototype {\em fractal landscape} that models the situation described in the introduction.

To define formally this function one can use the technique of Iterated Function Systems \cite{hut,dem} but in the present case we find it easier to proceed directly: write any number $x \in [0,1]$ via its base--three representation
\begin{equation}
\label{eq-hm0}
  x = \sum_{j=1}^\infty a_j 3^{-j},  \;\; a_j \in \{0,1,2 \}
\end{equation} 
and let $H^m$, for $m=1,\ldots$, be the subset of $[0,1]$ composed of all points such that the first one in their representation is found at the $m$-th position:
\begin{equation}
\label{eq-hm1}
   H^m = \left\{  x \mbox{  s.t. } a_j \neq 1 \mbox{  for } 1 \leq j < m, a_m=1  
   \right\}.
\end{equation}
For instance, $H^1 =(\frac{1}{3},\frac{2}{3})$, $H^2 =(\frac{1}{9},\frac{2}{9}) \cup (\frac{7}{9},\frac{8}{9})$. Clearly, these sets are the gaps of integer order $m$ in the usual construction of the Cantor set, obtained by deleting middle thirds. Equivalently, the set $H^m$ is the set of middle--third intervals discarded at level $m$ in the ternary Cantor set construction.  Therefore,
$H^m$ is a collection of open intervals, of length $3^{-m}$ and cardinality $2^{m-1}$. For $x \in [0,1]$ the discrete function $\mathbf{m}$ 
\begin{equation}
\label{eq-dis1}
   \mathbf{m}(x) = \left\{ \begin{array}{ll} m, &  x \in  H^m, m = 1,\ldots, \\
   \infty & \mbox{ otherwise.  }
   \end{array} \right.
\end{equation}
keeps track of the set $H^m$ to which it belongs. Clearly, 
for Lebesgue almost all $x \in \chi$ the order $\mathbf{m}(x)$ is finite, we set it to infinity in the residual set.

We stipulate that the physical function $f$ depends only on $\mathbf{m}(\cdot)$:
\[
f(x) = \tilde{f}(\mathbf{m}(x))  
\]
and it is therefore a piece-wise constant function, whose value is constant on every set $H^m$ (see Fig. \ref{disegno}).
The monotonically increasing function $\tilde{f}$ of the integer argument $m$ will be specified shortly. Whatever this choice,
the function $f$ is discrete and is characterized by a countable infinity of {\em plateaus} on the gaps of a Cantor set. On the Cantor set itself, we take $f$ and ${\bf m}$ to be infinity.

It is now time to specify the dynamics. We examine two different maps, both of which leave the Lebesgue measure $d\mu = dx$ invariant. The first is the hyperbolic dynamics generated by the asymmetric tent map: $\varphi_1(x) = x/p$, for $0\leq x < p$, and $\varphi_1(x) = (1-x)/(1-p)$ for $p\leq x \leq 1$, with the parameter $p$ in $(0,1)$. The second is the irrational rotation $\varphi_2 : x \rightarrow x + \omega \mbox{ mod }(1)$, with $\omega \in \R \setminus \Q$. As it is well known, this dynamics does not exhibit any decay of correlations.

On this basis, let us carry on the analysis sketched in Section \ref{sec-extr}.
Letting $m \in \bna$ and $h_m = \tilde{f}(m)$, we easily find that
\begin{equation}
\label{barf}
\bar{F}(h_m)= \mu \{ x \in \chi \mbox{  s.t. } f(x) \geq h_m \} = (2 \delta)^m.
\end{equation}
The value $\delta=1/3$ is the self--similarity ratio of the ternary Cantor set.
Notice that the function $\overline{F}$ is discontinuous and we are in the same situation discussed in \cite{nonsmooth}: we adopt here the subsequence approach, to consider a discrete set of values  $h_m$ and $\tau_n$ (see below). Also notice that $m$ is a label of the different level sets of $f$ and therefore of the intensity $f(x)$ associated to the phase-space point $x$. Before making explicit reference to the function $\tilde{f}$, let us derive the underlying extreme value law associated with this dynamical process, assuming heuristically independence of the dynamical variables, which can be justified under the strong mixing properties of the dynamical system $\varphi_1$.
We now write eq. (\ref{eq-2}) as
\begin{equation}
\label{eq-ca1}
   \tau_n = n \bar{F}(h_{m(n)}),
\end{equation}
where $m$ is a function of $n$: $m = m(n)$.
In order that $\tau_n \to \tau$, as in eq. (\ref{eq-2}), we need that
\begin{equation}
\label{eq-ca2}
  \log (\tau_n) =  m(n) \log 2 \delta + \log n \to \log \tau
\end{equation}
{\em i.e.}
\begin{equation}
\label{eq-ca3}
    m(n) = \frac{\log n/ \tau}{\log 1/2 \delta} + o(n).
\end{equation}
Let us henceforth drop the $o(n)$ term in $m(n)$ and define the inverse function $n = n(m)$ that yields the discrete sequence of values $\tau_n$ that can be realized, since the variable $m$ is discrete.
Then, we expect an extreme value law given by eq. (\ref{eq-3bis}), with $h_n(y)$ replaced by $h_{m(n)}$.

The threshold values $h_{m(n)}$ depend on the specific choice of the function $\tilde{f}$. Yet, the existence of an extreme value law depends on the discrete nature of this process and can be verified prior to such choice. In fact, let $\mathbf{X}_j = \mathbf{m} \circ \varphi^j$ and $\mathbf{M}_n=\max\{\mathbf{X}_0,\ldots,\mathbf{X}_{n-1}\}$,
That is, we consider the extreme value of the integer function $\mathbf{m}$ along a trajectory. Then, it is clear that
$
    \{ x \mbox{  s.t.  } M_n(x) \leq h_{m(n)} \} =
\{ x \mbox{  s.t.  } \mathbf{M}_n(x) \leq m(n) \}
$
and eq. (\ref{eq-3bis}), with $\theta=1$, is equivalent to
\begin{equation}
\label{eq-ca42}
   A(m,n) := \mu (\{ x \in \chi \mbox{  s.t.  } \mathbf{M}_n(x) \leq m \}) \rightarrow e^{-n (2 \delta)^m}.
\end{equation}

\begin{figure}
\centerline{\includegraphics[width=.6\textwidth, angle = -90]{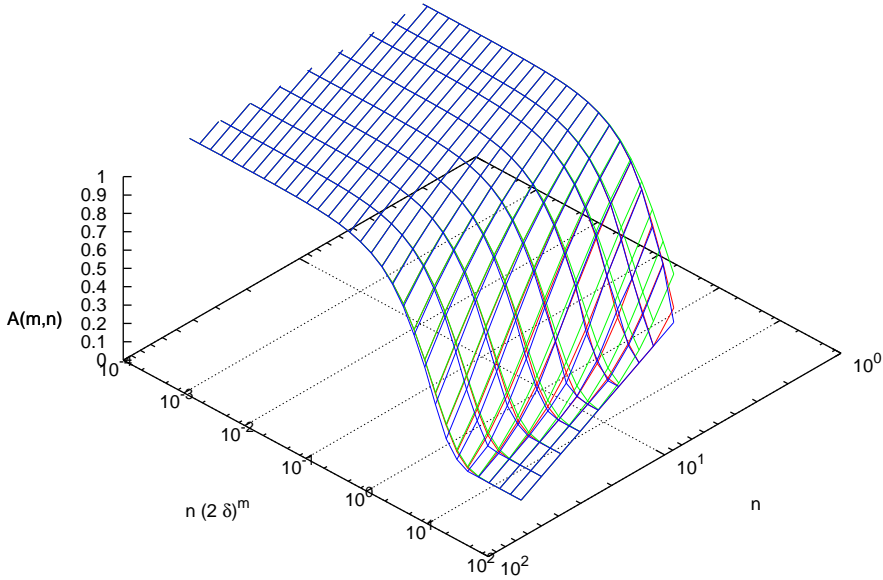}}
\caption{Experimental function $A(m,n)$ defined in eq. (\ref{eq-ca42}) for the dynamics $\varphi_1$ (red), the dynamics $\varphi_2$ (green). Lines in the grid connect raw data points at constant $n$ and at constant $m$: no surface interpolation is performed. 
The value of $m$ grows on lines at constant $m$ when moving towards the left of the figure, where $A(m,n)$ tends to one. Also reported is the graph of the function $e^{-n (2 \delta)^m}$ (blue). The three plots are almost coincident:  Figure \ref{fig3b.ps} draws their difference.}
\label{fig1c.ps}
\end{figure}

\begin{figure}
\centerline{\includegraphics[width=.6\textwidth, angle = -90]{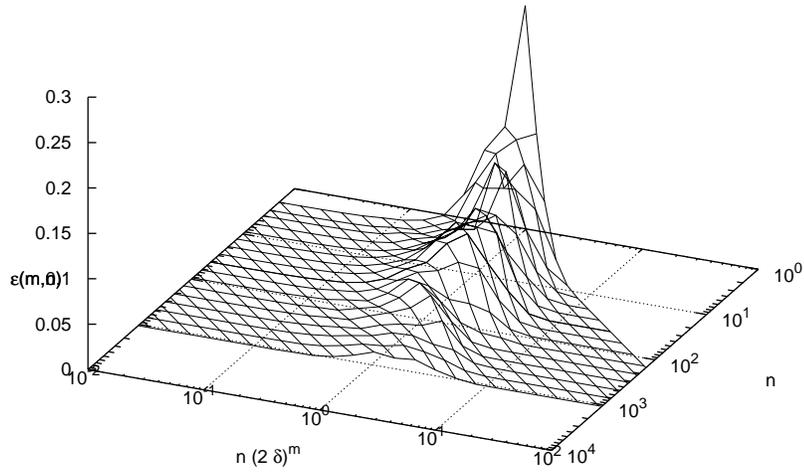}}
\caption{Difference $\varepsilon(m,n)$ between dynamical data $A(m,n)$ for $\varphi_2$ and the limit function $e^{-n (2 \delta)^m}$. As in Fig. \ref{fig1c.ps} the plot is rendered by lines connecting experimental data at constant $n$ and $m$. The line with $m=0$ is the first to appear to the right  of the figure.  }
\label{fig3b.ps}
\end{figure}

This dependence can be verified by plotting $A(m,n)$ versus $n$ and the combined variable $\tau_n=n (2 \delta)^m$: in Figure \ref{fig1c.ps} we observe convergence to the extreme value law, as $n$ grows.
Observe that {\em three} sets of data are reported in the figure and appear as almost undistinguishable---see the successive Figure \ref{fig3b.ps} to appreciate differences. The first is the function $e^{-\tau_n}=e^{-n (2 \delta)^m}$ (that appears as a pure negative exponential, with the choice of axes  of the figure); the second is the plot of the numerical data $A(m,n)$ for the dynamical system $\varphi_1$ described above; the third is analogous to the second, but under the dynamics $\varphi_2$. It is remarkable that also in this merely ergodic case the extreme value law is obtained in the limit. Also observe that the extremal index takes the value one, in both cases.

{\em In lieu} of a formal justification of these observations, which we defer to future works, we conjecture that the r\^ole of fast decay of correlations in assuring the dependence conditions \cite{jorge1} is played here by the complexity of the intensity function, a fact that could also explain the absence of clustering of extreme events.

These results are better appreciated in Figure \ref{fig3b.ps}, where we plot the difference $\varepsilon(m,n):= |A(m,n)-e^{-n (2 \delta)^m}|$ for the case of the map $\varphi_2$, which is larger than that observed for the  case $\varphi_1$. As $n$ grows $\varepsilon(m,n)$ tends to zero for all values of $m$. In this paper, we are not concerned with the rate of convergence to the limit distribution that could nonetheless be derived from these data.
In the conventional framework, convergence rates have been studied in \cite{speed}.

\begin{figure}
\centerline{\includegraphics[width=.6\textwidth, angle = -90]{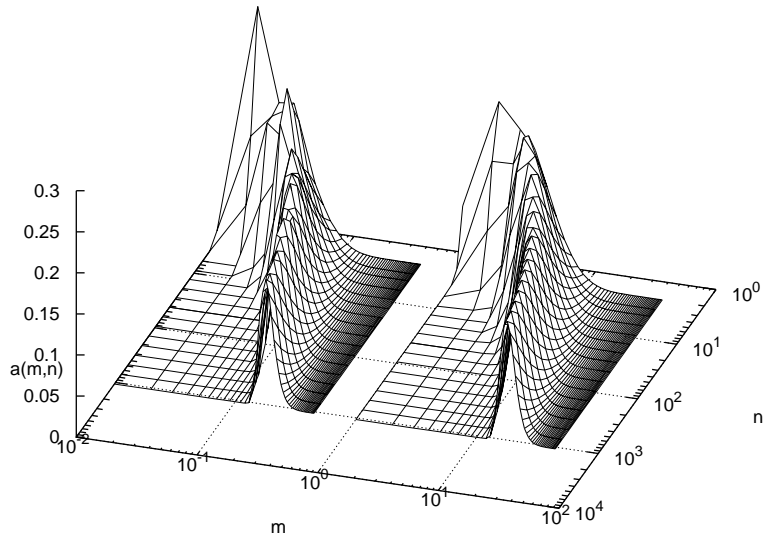}}
\caption{Empirical densities $a(m,n)$ (defined in the text) for the dynamics $\varphi_1$ (right plot) and $\varphi_2$ (left plot, rendered by the shift $m \rightarrow 10^{-2} m$).}
\label{fig4c.ps}
\end{figure}

Finally, in Figure \ref{fig4c.ps} we also plot side by side the empirical densities $a(m,n) := \mu (\{ x \in \chi \mbox{  s.t.  } \mathbf{M}_n(x) = m \})$ for the tent map $\varphi_1$ and the irrational rotation $\varphi_2$. The differences between the two fade away as $n$ grows and as they both tend to the limit extreme value law.

We can put the above results in the usual casting of extreme laws by letting $a=-\log(2 \delta) >0$ and thus re-writing eq. (\ref{eq-ca42}) as
\begin{equation}
\label{eq-ca50}
   \mu (\{ x \in \chi \mbox{  s.t.  } a [\mathbf{M}_n(x) - \log(n)/a ] \leq y \}) \rightarrow e^{-e^{-y}},
\end{equation}
where the variable $y$ is restricted to take discrete values, multiples of $a$.

We can now return to the problem of specifying the physical function $f(x)$.
As typical in dynamical extreme value theory, we can choose different forms of the function $\tilde{f}$. Let us follow the standard classification given {\em e.g.} in ref. \cite{extremal}, eqs. (3.3) - (3.5); for reasons that will become clear in the following section, functions $\tilde{f}_i(m)$ corresponding to the examples $g_i(\zeta)$ given in \cite{extremal} can be obtained by writing $\tilde{f}_i(-\log \zeta)= g_i(\zeta)$ and making the substitution $m = -\log \zeta$. This yields
\begin{eqnarray}
\label{eq-trecas1a}
\tilde{f}_1(m) = m \\
\label{eq-trecas1b}
\tilde{f}_2(m) = e^{m/\beta} \\
\label{eq-trecas1c}
\tilde{f}_3(m) = D - e^{-m/\gamma}, \quad  D>0, \quad \gamma>0.
\end{eqnarray}

Case 1 has been just studied: in the other cases it is clear that the validity of the extreme value law for the underlying discrete process $\mathbf{X}_j$ entails that of the derived processes $X_j$. By straightforward calculation, since $\mathbf{M}_n \leq m \Leftrightarrow M_n \leq \tilde{f}_i(m)$ in eq. (\ref{eq-ca42}), we obtain that, letting
$\eta = -\beta \log(2 \delta) = a \beta > 0$ in case $2$ and $\zeta = a \gamma > 0$ in case $3$,
the extreme value law
\begin{equation}
\label{eq-trecas3}
\mu(\{ x \in \chi \mbox{  s.t.  } a_j(n) [ M_n(x) -b_j(n)] \leq y\}) \rightarrow e^{-\theta \tau(j;y)}
\end{equation}
holds with
\begin{eqnarray}
\label{eq-trecas1}
a_1(n) = a, \; b_1(n) = \log(n) /a,\; \tau(1;y) = e^{-y};   \\
a_2(n) = n^{-1/\eta}, \; b_2(n) = 0, \; \tau(2;y) = y^{-\eta}; \\
a_3(n) = n^{1/\zeta} , \; b_3(n) = D , \; \tau(3;y) = (-y)^{\zeta};
\end{eqnarray}
that orderly correspond to the Gumbel, Fr\'echet and Weibull forms, respectively.

Let us now take a moment to comment on the second case.
Observe that this construction accounts for the Gutenberg--Richter law, according to which the statistical frequency of events more intense than a threshold $y$ is a power law in $y$. In fact, by choosing $\tilde{f} = \tilde{f}_2$ we obtain $h_m = e^{m/\beta} = y$, which implies that
\[
\bar{F}(y)= \mu \{ x \in \chi \mbox{  s.t. } f(x) \geq y \} = y^{\beta \log( 2 \delta)} = y^{-\eta}.
\]
Ergodicity of the dynamical systems with maps $\varphi_1$ and $\varphi_2$ and Lebesgue measure clearly implies that the phase-space average above is equal to a time average.
The exponent $\eta = -\beta \log (2 \delta)$ is experimentally measured in a wide set of natural phenomena (earthquakes, avalanches, floods, {\em et cetera}).
We can so conclude that the simple model based on the Cantor Ladder well describes events with Gutenberg--Richter statistics and yields a Frechet law for extrema parameterized by the Gutenberg--Richter exponent. Relation between the two laws is well known in the natural sciences: for an example in geophysics, see
\cite{lavenda}.

\section{Rare events: Minkowski dimension and content}
\label{sec-discantor}

In the previous section we have described a simple dynamical model for extreme values that differs from standard usage in the literature in a significant respect: the intensity function $f$ is {\em not} related to the distance of the motion $x_j = \varphi^j(x_0)$ from a point $x^*$ (or to the measure of a ball centered at $x^*$), but it is a global function with a hierarchical structure. In this section we show that such model can be seen as an extension of the standard approach. We first show the analogy and we then develop the new theory.

Consider again the {\em gap order function}
$\mathbf{m}(x)$, that yields $f(x) = \tilde{f}(\mathbf{m}(x))$. Recall that $\mathbf{m}(x) = m$ if and only if $x$ belongs to $H^m$, {\em i.e.} to a gap of order $m$ in the hierarchical construction of the Cantor set. Since the Cantor sets that we are considering are of null Lebesgue measure, almost all points $x_j$ of a typical trajectory belong to a gap. Then, the closest point of the Cantor set $K$ to $x_j$ is one of the two gap extrema. Define now the distance of a point $x \in \chi$ to a Cantor set $K$ as
\begin{equation}
\label{eq-cant1}
    d(x,K) = \min \{ d(x,s), s \in K \},
\end{equation}
where compactness of $K$ permits to substitute the required $\inf$ with a  minimum. The previous consideration permits us to compute  the average value, with respect to the Lebesgue invariant measure, of the {\em logarithm} of the distance $d(x,K)$ when $x$ belongs to a gap $G$ of length/Lebesgue measure $l$:
\begin{equation}
\label{eq-cant2}
    \frac{1}{l} \int_G \log (d(x,K)) d\mu(x) = \log (l) - c,
\end{equation}
where $c=1+\log 2$.
In our simple construction, when $G \in H^m$, the diameter of the gap is $l=\delta^m$ and this proves that $\mathbf{m}(x)$ {is proportional to the logarithm of the average distance of a point in a gap, from the Cantor set $K$.}

Let us introduce the following generalization of the standard setting: rather than defining $f(x)$ via the distance of $x_j$ from a finite set of points, replace the latter by a generic compact set $K$. The fundamental definition now becomes
\begin{equation}
\label{eq-cant4}
    f(x) = \tilde{f}_j (-\log (d(x,K))),
\end{equation}
in which $\tilde{f}_j$ is any of the three forms in eqs. (\ref{eq-trecas1a})--(\ref{eq-trecas1c}).
In force of this analogy, we {\em conjecture that the extreme value laws obtained in the previous section for the Cantor Ladder model should correspond orderly to the choices of $\tilde{f}$ in eq. (\ref{eq-cant4}), when $K$ is the ternary Cantor set.}

Without loss of generality, let us henceforth concentrate on the case $\tilde{f}_1$.
The main object in the theory is the set where $f$ exceeds a threshold $h$:
\begin{equation}
\label{eq-cant5}
     \{ x \in \chi \mbox{  s.t. } f(x) \geq h \} = \{ x \in \chi \mbox{  s.t. } d(x,K) \leq e^{-h} \} =
     N_{e^{-h}}(K).
\end{equation}
In the new setting, this becomes the {\em neighborhood} of radius $e^{-h}$ of the set $K$, a well known object in fractal geometry. Clearly, in the case $K=\{x^*\}$ this reduces to the frequently investigated ball around the point $x^*$, indicated by $B_{e^{-h}}(x^*)$.
Taking the measure of the sets in eq. (\ref{eq-cant5}) we obtain
\begin{equation}
\label{eq-expoball}
     \overline{F}(h) =  \mu (N_{e^{-h}}(K)).
\end{equation}
Let us now consider again the case when the invariant measure $\mu$ in the above equation is the Lebesgue measure and let
\begin{equation}
\label{eq-cant6}
{d}_M (K) = \lim_{\varepsilon \to 0}
     \frac{ \log \mu (N_{\varepsilon}(K))}{\log \varepsilon}.
\end{equation}
The {\em Minkowski dimension} of a set $K$ is defined as the difference $D-d_M(K)$, where $D$ is the Euclidean dimension of the ambient space \cite{lapi}. It may, or it may not, coincide with the Hausdorff dimension: it does, for the ternary Cantor set; we will see an instance of the opposite case in section \ref{sec-usuenne}. Upper and lower dimensions are defined in the standard way when the limit in eq. (\ref{eq-cant6}) does not exist.
Furthermore, a set is defined to be {\em Minkowski measurable} if the limit
\begin{equation}
\lim_{\varepsilon \rightarrow 0}
 \mu (N_{\varepsilon}(K))  \varepsilon^{-d_M (K)} = A_{K}
 \label{eq-dmin1c}
\end{equation}
exists and is neither zero nor infinity. The constant $A_{K}$ is called the {\em Minkowski content} of the set $K$.

If we now follow standard derivations in dynamical extreme value theory we find that the {\em local dimension} at the point $x^*$ is now replaced by $d_M(K)$. Recall that this is the case when the invariant measure of the dynamics, $\mu$, is the Lebesgue measure. We will generalize this situation in the following section.
Let us henceforth suppose that an asymptotic law holds, of the kind
\begin{equation}
\label{eq-cant61}
 \mu (N_{\varepsilon}(K)) \sim A_K \varepsilon^{d_M (K)},
\end{equation}
where $A_K$ is now the {average} Minkowski content of the set $K$\footnote{It is to be remarked that cases like the ternary Cantor set studied in this paper belong to what are called IFS of lattice type. For these IFS the limit in eq. (\ref{eq-dmin1c}) exists only in a suitable Ces\'aro average, due to characteristic log-periodic oscillations \cite{daniel,guzzet}. It is then called the average Minkowski content \cite{lapi}.}. Then,
\begin{equation}
\label{eq-cant62}
 \overline{F}(h) \sim A_K e^{- h d_M (K)}
\end{equation}
and eq. (\ref{eq-2}) implies that
\begin{equation}
\label{eq-cant62z}
 \tau_n = n A_K  e^{-h d_M (K)},
\end{equation}
\begin{equation}
\label{eq-cant62b}
 h_n = \frac{1}{d_M (K)} \log ( \frac{A_K n}{\tau} ).
\end{equation}
These equations give the set of thresholds to be used in eq. (\ref{eq-3bis}):
we now test numerically the extreme value law.

\begin{figure}
\centerline{\includegraphics[width=.4\textwidth, angle = -90]{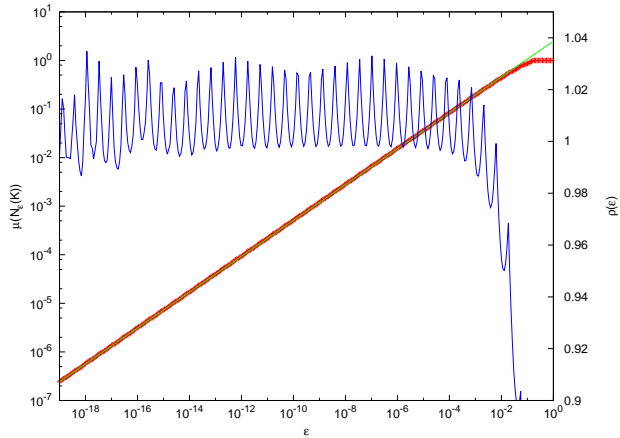}}
\caption{Measure of the $\varepsilon$ neighborhood of the Cantor set $K$ versus $\varepsilon$ (red crosses) and asymptotic function $A_K \varepsilon^{d_M (K)}$ (green line), with $d_M(K) = 1 -\log(2)/\log(3)$ and $A_K = 5/2$. The numerical sample has been computed via the irrational rotation dynamics $\varphi_2$. Also plotted is the ratio $\rho(\varepsilon)$ of  numerical sample and theoretical power--law (blue line, right vertical scale). }
\label{fig5b.ps}
\end{figure}

Let $K$ be the ternary Cantor set previously considered. In Figure \ref{fig5b.ps} we plot, in double logarithmic scale, the function $\mu (N_{\varepsilon}(K))$ versus $\varepsilon$, together with the asymptotic estimate of eq. (\ref{eq-cant61}). 
Two remarks are in order at this point. The first is that  parameters in eq. (\ref{eq-cant61}) can be computed exactly: $d_M(K) = 1 -\log(2)/\log(3)$ and $A_K = 5/2$. The second is that {\em log--periodic} oscillations in $\varepsilon$ are super--imposed to the scaling behavior (and the same happens for the statistics $A(h,n)$). Nonetheless, these oscillations, albeit theoretically important, are of inferior practical significance than the leading behavior.

\begin{figure}
\centerline{\includegraphics[width=.6\textwidth, angle = -90]{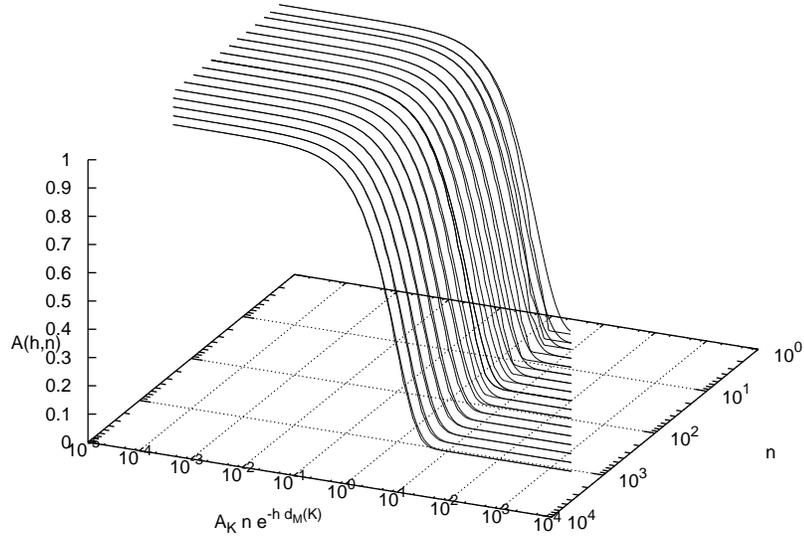}}
\caption{Experimental extreme value function $A(h,n)=\mu (\{M_n < h \})$ for the Cantor set $K$ and irrational rotation dynamics $\varphi_2$, versus $n$ and $\tau_n = A_K n e^{-h d_M(K)}$, together with $e^{-\tau_n}$. Already at $n=1,000$ the two curves are almost undistinguishable. Differences are plotted in Fig. \ref{newfig6a2.ps}.}
\label{fig6b.ps}
\end{figure}

As in the previous section, we can then plot $A(h,n) = \mu (\{x \in \chi \mbox{  s.t. } M_n(x) < h \})$ versus $n$ and $\tau_n$, to verify the asymptotic law
\begin{equation}
\label{eq-cant7}
\mu (\{x \in \chi \mbox{  s.t. } M_n(x) < h \}) \rightarrow \exp \{ - \theta n A_K e^{-h d_M(K)}\},
\end{equation}
in which convergence is again to be understood in leading order, with super--imposed oscillations.
This is done in Figure \ref{fig6b.ps} in the case of irrational rotations, $\varphi_2$. Convergence to the asymptotic law is evident. As was to be expected from the heuristic argument of the previous section, the extremal index $\theta$ is again one. Finally, we can resort to normalizing sequences, to show that we are in the presence of a Gumbel type (suitably translated results hold for the cases $\tilde{f}_2$ and $\tilde{f}_3$) extreme value law :
\begin{equation}
\label{eq-cant71}
\mu (\{x \in \chi \mbox{  s.t. }  d_M(K) M_n(x) - \log (A_K n) < y \}) \rightarrow \exp \{ -  e^{-y}\},
\end{equation}
in which the r\^ole of the Minkowski dimension {\em and} of the asymptotic constant $A_K$, the Minkowski content, is evident. This is a remarkable case in which these geometrical quantities take on a dynamical meaning.

\begin{figure}
\centerline{\includegraphics[width=.6\textwidth, angle = -90]{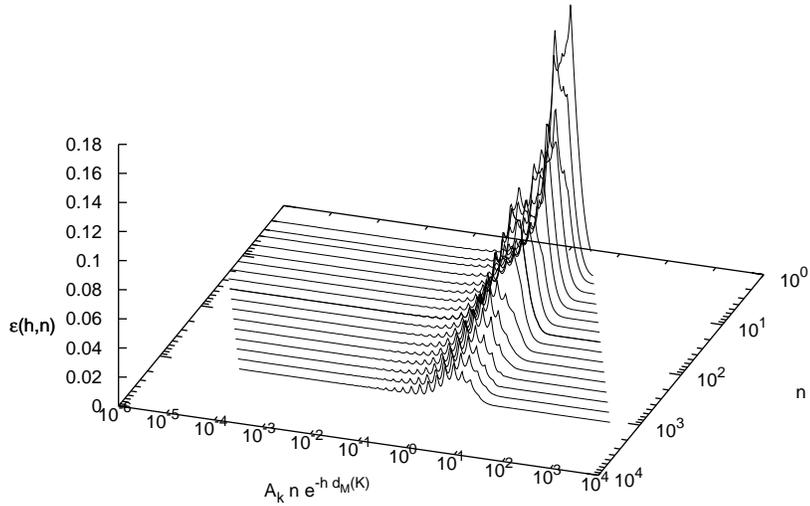}}
\caption{Absolute difference $\varepsilon(h,n)=|A(h,n)-e^{-\tau_n}|$ for the Cantor set $K$ and irrational rotation dynamics $\varphi_2$, versus $n$ and $\tau_n = A_K n e^{-h d_M(K)}$.}
\label{newfig6a2.ps}
\end{figure}

\section{Generalized Minkowski dimension and content}
\label{sec-genmink}

So far, we have considered dynamical systems whose invariant measure is the Lebesgue measure. It is now interesting to consider measures significantly different from the latter: in particular, those that are singular continuous with respect to it.
For sake of simplicity and with no loss of generality,
we still consider the case of the previous section, when  $f(x) = -\log d(x,K)$. As a consequence, recalling eq. (\ref{eq-cant5}), the function $\overline{F}(h) =  \mu (N_{e^{-h}}(K))$ may now depend on the measure $\mu$ in a non--trivial fashion.

To begin with, suppose that a relation like eq. (\ref{eq-cant61}) holds,
\begin{equation}
 \mu (N_{\varepsilon}(K)) \sim A_{K,\mu} \varepsilon^{d_M (K,\mu)},
 \label{eq-dmin1}
\end{equation}
where $d_M(K,\mu)$ and $A_{K,\mu}$ are now functions of the relation between the set $K$ and the measure $\mu$. We therefore define the {\em generalized Minkowski dimension} as the difference between the Euclidean dimension of the embedding space and $d_M(K,\mu)$. Likewise, the constant $A_{K,\mu}$ defines the {\em generalized Minkowski content}, as in eq. (\ref{eq-dmin1c}).
Under these circumstances, eq. (\ref{eq-cant62}) still holds, with the obvious replacement of symbols, and we may expect that the previous conclusions on the extreme events statistics, eqs. (\ref{eq-cant62z})--(\ref{eq-cant7}), also apply.

Let us now introduce an important dynamical system, defined by the map $\varphi_3$:
\begin{equation}
\label{eq-mink1}
    \varphi_3(x) = \left\{ \begin{array}{ll}
    \frac{x}{1-x} & 0 \leq x \leq \frac{1}{2},\\
    2 - \frac{1}{x} & \frac{1}{2} < x \leq 1.
    \end{array}
    \right.
\end{equation}
This map is a sort of non-linear version of the $x \rightarrow 2x \mbox{ mod } 1$ map, this time composed of M\"obius transformations. It has a singular continuous invariant measure $\mu$  that appears in many problems of number theoretical origin \cite{guzzi}. In dynamical systems, it has been studied for instance in \cite{myprl}, where a graph of its integrated density, $Q(x)=\int^x d\mu(s)$ can be found. The function $Q$ is also known as the Minkowski's question--mark function \cite{minkoqm}. It satisfies the relations
\begin{equation}
\label{eq-mink2}
       \begin{array}{lcl}
    Q(\frac{x}{1+x}) & = & \frac{1}{2} Q(x) ,\\
         Q(\frac{1}{2-x}) & = & \frac{1}{2} [Q(x)+1],
    \end{array}
\end{equation}
that equivalently serve to define it via a M\"obius Iterated Function System \cite{myprl}.
Among its remarkable properties one notices the fact that the derivative of $Q$ exists and is null almost everywhere, but never on a full interval. In addition, $Q$ has a non--analytical behavior at all rational points, as we will see below.

In Figure \ref{figc34f15} we plot the measure of the $\varepsilon$ neighborhood of the ternary Cantor set $K$ versus $\varepsilon$, and the asymptotic law (\ref{eq-dmin1}), in which $d_M(K,\mu) = 1 - \log(2) / \log (3)$ and $A_{K,\mu}=1$. That is to say, for this measure of high theoretical importance, the ternary Cantor set has Minkowski dimension $\log(2)/\log(3)$ (as in the Lebesgue case) but average Minkowski content one!

\begin{figure}
\centerline{\includegraphics[width=.4\textwidth, angle = -90]{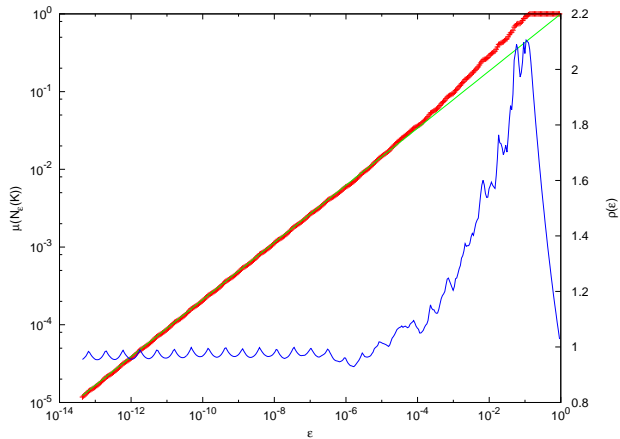}}
\caption{Measure of the $\varepsilon$ neighborhood of the Cantor set $K$ versus $\varepsilon$ (red points) and asymptotic function $A_{K,\mu} \varepsilon^{d_M (K,\mu)}$ (green line), with $d_M(K,\mu) = 1 -\log(2)/\log(3)$ and $A_{K,\mu} = 1$. The numerical sample has been computed via the M\"obius dynamics $\varphi_3$ with invariant measure $d\mu=dQ$. Also plotted is the ratio $\rho(\varepsilon)$ of  numerical sample and theoretical power--law (blue line, right vertical scale).}
\label{figc34f15}
\end{figure}

\begin{figure}
\centerline{\includegraphics[width=.6\textwidth, angle = -90]{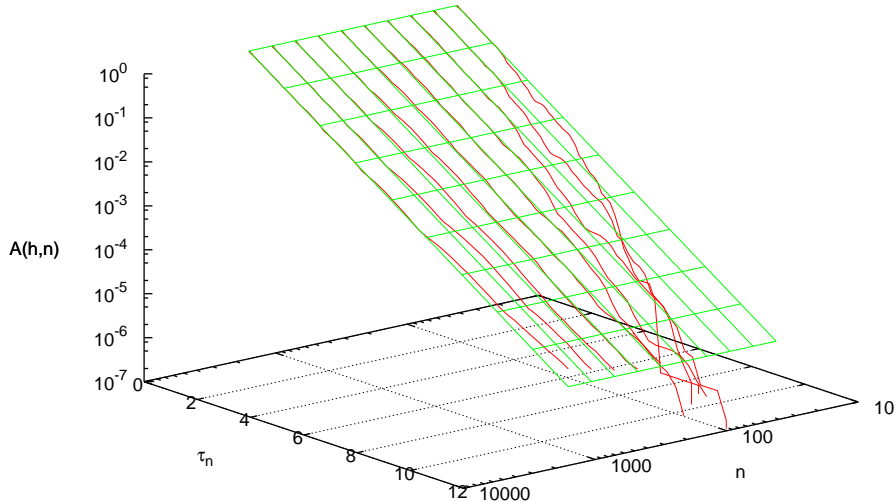}}
\caption{Experimental extreme value function $A(h,n)=\mu (\{x \in \chi \mbox{  s.t. } M_n(x) < h_n(\tau) \})$ (red lines) for the ternary Cantor set $K$, M\"obius dynamics $\varphi_3$, and invariant Minkowski measure $d\mu = dQ$, versus $n$ and $\tau_n$ (eq. (\ref{eq-cant62z})), together with $e^{-\tau_n}$ (green lines).}
\label{figc34-13}
\end{figure}

Let us now examine the dynamical motions induced by the map $\varphi_3$ and compute the extreme value statistics. We obtain the data reported in Fig. \ref{figc34-13} where, at difference with the other figures, we have employed a linear scale for $\tau_n$ and a logarithmic scale for $A(h,n)$. Convergence to the extreme value law is observed over many orders of magnitude, until fluctuations due to finite sampling emerge (6,400,000 samples of the variable $M_n(x_j)$ have been computed for each value of $n$).

The results of this section can be explained heuristically by saying that the motion of the dynamical system $(\chi,\varphi_3,Q)$ samples space around the Cantor set $K$ almost like the Lebesgue measure, the difference emerging in the Minkowski content while not in the Minkowski dimension. This is certainly an important observation, that adds to the long list of non--trivial behaviors of this dynamical system and calls for rigorous study.

\section{Extremely rare events: non--standard Minkowski scaling}
\label{sec-extrare}

Extreme value laws are linked to the statistics of rare events.
In fact, consider for instance the dynamics studied in \cite{towards}, where a selected point is chosen on the H\'enon attractor, at which the intensity of the simulated phenomenon is maximal. Rare events occur when the motion enters balls around this location. Typically, the measure of the ball scales in a power--law fashion with its radius and the local dimension of the invariant measure exists as a positive number. We have extended this set--up to neighborhoods of a set $K$, via the generalized Minkowski dimension and content. When these latter exist, results like those of the previous sections are to be expected.

Yet, there are cases where the interplay between the exploring dynamical system and the compact set $K$ of singularities of the intensity function (embedded, or not, in its attractor) is more intricate than what observed at the end of the previous section. In these cases, different scaling properties than the one given in eq. (\ref{eq-dmin1}) may hold, implying important differences in the extreme value distribution. In this section we describe one of these cases that can be interpreted as describing the statistics of {\em extremely rare events}: in fact, in eq. (\ref{eq-expoball}) not only the radius of the set scales exponentially with the intensity $h$ of the phenomenon, but also its measure scales much more rapidly than power--law in the radius.

Consider again the dynamical system $(\chi,\varphi_3,Q)$ defined in the previous section and the set $K=\{0\}$.
We know that, for small $x$, $Q(x) \sim 2 e^{-\frac{\log 2}{x}}$, as can be shown using eq. (\ref{eq-mink2}). Since $Q(x)=\mu(B_x(0))=\mu(N_x(K))$, we are in the presence of a neighborhood of a set $K$ whose measure does not scale as
eq. (\ref{eq-dmin1}). Similarly, when $x=1/2$, one computes $\mu(B_\varepsilon(1/2)) = Q(\frac{1}{2}+\varepsilon) - Q (\frac{1}{2}-\varepsilon) \sim \sqrt{2} \exp \{ - \frac{\log 2}{4 \varepsilon}\}$. More generally, still using using eq. (\ref{eq-mink2}), it is possible to derive an asymptotic formula for the measure of all balls centered at points of the kind $\frac{1}{k}$ (valid for $k>>1$ and small $\varepsilon$):
\begin{equation}
\label{eq-mink3}
    \mu(B_\varepsilon(\frac{1}{k})) \sim 2^{3-k-1/k} \; e^{-\frac{\log 2}{k^2 \varepsilon}}.
\end{equation}
We are therefore in the presence of a dynamical system and a family of sets $K$,  the singletons $\{\frac{1}{k}\}$, for which an asymptotic behavior different from eq. (\ref{eq-dmin1}) holds: precisely,
\begin{equation}
\label{eq-mink4}
    \mu(N_\varepsilon(K))  \sim B_{K,\mu} \; e^{-D_{K,\mu} \varepsilon^{-q_{K,\mu}}}.
\end{equation}
We would like to call the positive parameters $D_{K,\mu}$, $B_{K,\mu}$ and $q_{K,\mu}$ the {\em non--standard, generalized Minkowski constants}. Of course, they no more possess the properties of a dimension. In eq. (\ref{eq-mink3}) one has $q_{K,\mu}=1$: the wider generality will be needed in the next section.
Consequences to the extreme value statistics are then important:
\begin{equation}
\label{eq-mink41}
    \overline{F}(h)  \sim B_{K,\mu} \; e^{-D_{K,\mu} e^{hq_{K,\mu}}};
\end{equation}
therefore
\begin{equation}
\label{eq-mink42}
   \tau_n =   n B_{K,\mu} \; e^{-D_{K,\mu} e^{hq_{K,\mu}}}
\end{equation}
and
\begin{equation}
\label{eq-mink43}
   h_n(\tau) =  - \frac{1}{q}
   \log \left[  \log{ (\frac{ n B_{K,\mu}}{\tau})^{\frac{1}{D_{K,\mu}}} }    \right].
\end{equation}

Observe the double exponential in the $h$ dependence of $\tau_n$, or the double logarithmic dependence of $h_n$ on $\tau$. They imply a sharp decay of the distribution function of the extremal variable $M_n$, whose typical value increases with the logarithm of the logarithm of the time interval $n$. In most realistic applications, this increase could be practically undetectable: it is a consequence of the {\em extremely rare event} $N_\varepsilon(K)$.

\begin{figure}
\centerline{\includegraphics[width=.6\textwidth, angle = -90]{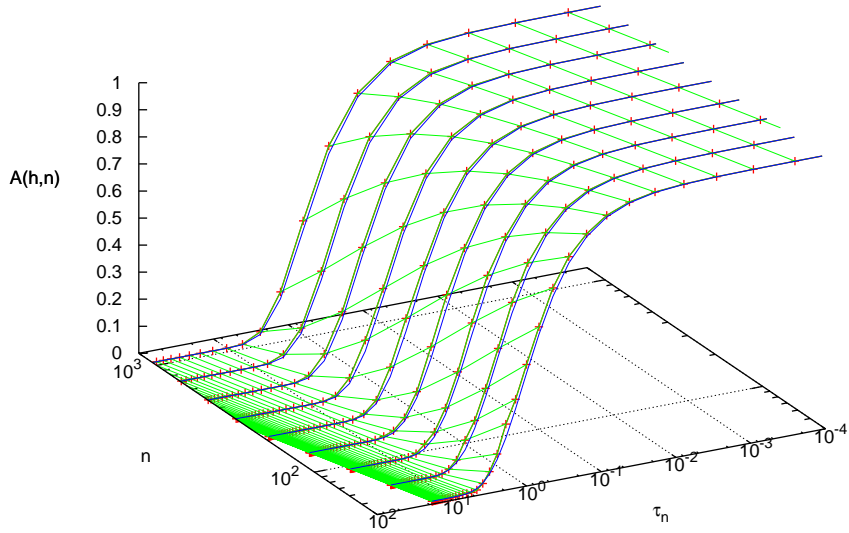}}
\caption{Experimental extreme value function $A(h,n)=\mu (\{x \in \chi \mbox{  s.t. } M_n(x) < h_n(\tau) \})$ (red crosses) for $K=\{1/4\}$, M\"obius dynamics $\varphi_3$ and invariant Minkowski measure $d\mu = dQ$, versus $n$ and $\tau_n$ (eq. (\ref{eq-mink42}), together with $e^{-\tau_n}$ (blue lines). Experimental points corresponding to the same value of $h$ are joined by green lines. }
\label{figp25}
\end{figure}

\begin{figure}
\centerline{\includegraphics[width=.6\textwidth, angle = -90]{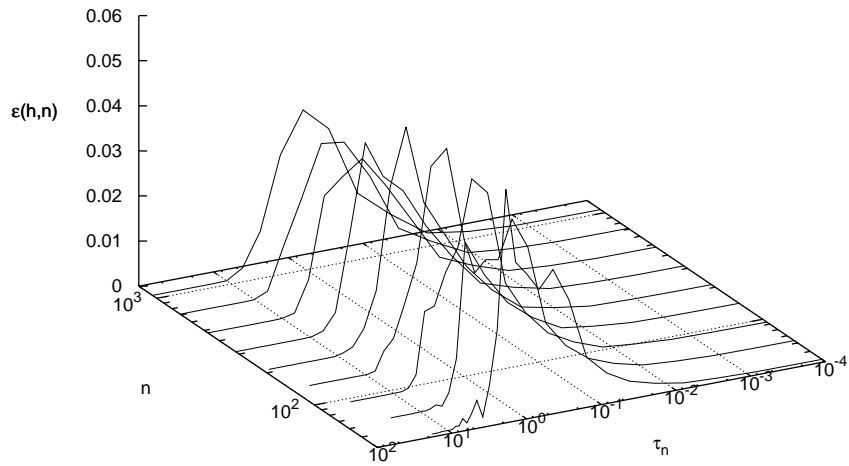}}
\caption{From Figure \ref{figp25}: $\varepsilon(h,n)=|A(h,n)- e^{-\tau_n}|$ for the M\"obius dynamics $\varphi_3$ and invariant Minkowski measure $d \mu = dQ$ versus $n$ and $\tau_n$. }
\label{figp25a}
\end{figure}

We can make use of the power of numerical experiments (as opposed to real experiments) to verify the previous conclusions: let us take $K=\{\frac{1}{k}\}$ and let us choose $k=4$.
As before, validity of eq. (\ref{eq-3bis}) can be verified by plotting the experimental extreme value distribution $A(n,h) = \mu \{x \in \chi \mbox{  s.t. } M_n(x) \leq h \}$ versus $n$ and $\tau_n$: this is done in Figure \ref{figp25}, with theoretical parameters $B_{K,\mu}=2^{-5/4}$, $D_{K,\mu}=\log(2)/16$ and $q_{K,\mu}=1$. The extremal index $\theta$ is one and the limiting curve (blue) is a pure exponential $e^{-\tau_n}$. In the same figure, we join by green lines the values corresponding to the same value of $h$. Difference between the experimental data and the theoretical limit curve are plotted in Fig. \ref{figp25a}.

Observe that the functional dependence in eq. (\ref{eq-mink43}) renders difficult to bring eq. (\ref{eq-3bis}) in one of the three standard forms and we are forced to parameterize the threshold $h_n$ as a function of $\tau$. Otherwise, allowing for a non--linear scaling, we can obtain a Gumbel law for the {\em exponential} of the extreme values:
\begin{equation}
\label{eq-mink44b}
  \mu\{x \in \chi \mbox{  s.t. } D_{K,\mu} e^{-q_{K,\mu} M_n(x)} -  \log (nB_{K,\mu}) \geq y \}  \rightarrow e^{-e^y}.
\end{equation}

\section{Extremely rare events with non--unit marginal index}
\label{sec-usuenne}

As an exception to our approach, in the previous section we have considered a family of intensity functions with singularities concentrated at a single point. In this section we return to wider generality, by considering a set $K$ of infinite cardinality, compact, although not Cantor, which provides an interesting family of extremely rare events.
This set is a renown geometrical example:
\begin{equation}
\label{eq-usenne}
    K = \{0\} \cup \{ \frac{1}{k}, k \in \N \}.
\end{equation}
As it is well known, its Hausdorff dimension is null (it is a countable set of points) while its Minkowski dimension is one half, as it can be easily shown.

Let us first compute its generalized Minkowski constants, when the measure $\mu$ is the derivative of the Minkowski question--mark function $Q(x)$. We have to compute the $\mu$--measure of the $\varepsilon$--neighborhood of the set $K$.
Observe that $N_\varepsilon(K)$ can be written as
\begin{equation}
\label{eq-mink60}
    N_\varepsilon(K) = [0,\delta_\varepsilon] \cup \bigcup_{k=k_\varepsilon}^\infty
    B_\varepsilon(\frac{1}{k}),
\end{equation}
where $\delta_\varepsilon \sim \sqrt{\varepsilon/2}$, $k_\varepsilon \sim \sqrt{2/\varepsilon}$. We can now use the results of the previous section to compute the measure of $N_\varepsilon(K)$: setting $\lambda := \log(2)$, we obtain
\begin{equation}
\label{eq-mink61}
    \mu(N_\varepsilon(K)) \sim    2 e^{-\lambda / \delta_\varepsilon} + \sum_{k=2}^{k_\varepsilon}  2^{3-k-1/k} \; e^{-\frac{\lambda}{k^2 \varepsilon}}.
\end{equation}
Consider now the summation at r.h.s.:
\begin{equation}
\label{eq-mink62}
     S_\varepsilon := 2^3 \; \sum_{k=2}^{k_\varepsilon}  e^{-\lambda[k+\frac{1}{k}+ \frac{1}{k^2 \varepsilon}]}.
\end{equation}
Its behavior can be evaluated by a saddle point technique, with the result that
$S_\varepsilon \sim B_{K,\mu} e^{-D_{K,\mu} \varepsilon^{-1/3}}$, where $B_{K,\mu}$ and $D_{K,\mu}$ are Minkowski constants (that we do not attempt to compute explicitly here), so that
\begin{equation}
\label{eq-mink62b}
    \mu(N_\varepsilon(K)) \sim    2 e^{-\lambda \sqrt{2/\varepsilon}} + B_{K,\mu} e^{-D_{K,\mu} \varepsilon^{-1/3}}.
\end{equation}
Observe that we have two different non--standard behaviors appearing in this equation, both of which with a non--unit exponent (as different from the previous section). Evidently, the second term is leading, when $\varepsilon$ tends to zero: we are therefore in the conditions of eq. (\ref{eq-mink4}), with the explicit determination $q_{K,\mu}=1/3$, that constitutes still another generalization brought along by the theory.

\begin{figure}
\centerline{\includegraphics[width=.6\textwidth, angle = -90]{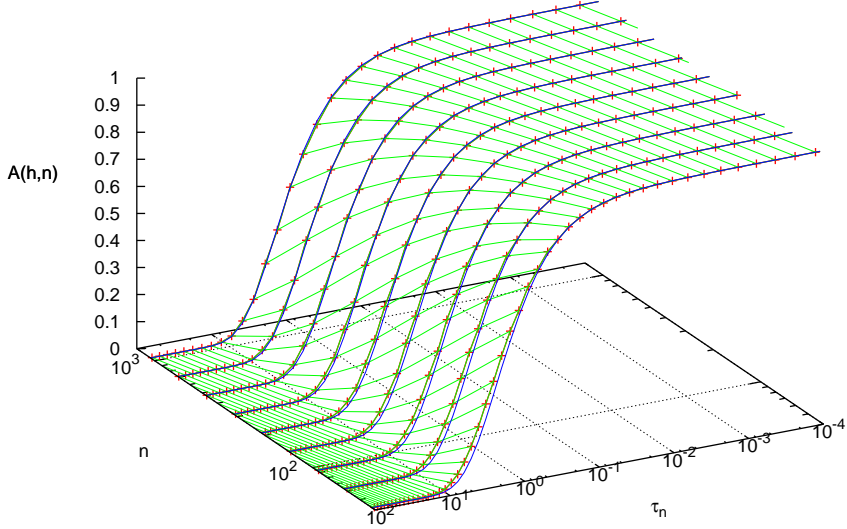}}
\caption{Experimental extreme value function $A(h,n)=\mu (\{M_n < h_n(\tau) \})$ (red crosses) for $K=\{0\} \cup \{ \frac{1}{k}, k \in N \}$, M\"obius dynamics $\varphi_3$ and invariant Minkowski measure $d\mu=dQ$, versus $n$ and $\tau_n$ in eq. (\ref{eq-mink42}), together with $e^{-\theta \tau_n}$ (blue lines). Experimental points corresponding to the same value of $h$ are joined by green lines. Experimental extremal index is $\theta \simeq 0.47$. Non-standard Minkowski constants take the values $B_{K,\mu}\simeq 24.61$ and $D_{K,\mu} \simeq 1.26$. The theoretical exponent $q_{K,\mu}=1/3$ has been used.}
\label{figp25b}
\end{figure}

\begin{figure}
\centerline{\includegraphics[width=.6\textwidth, angle = -90]{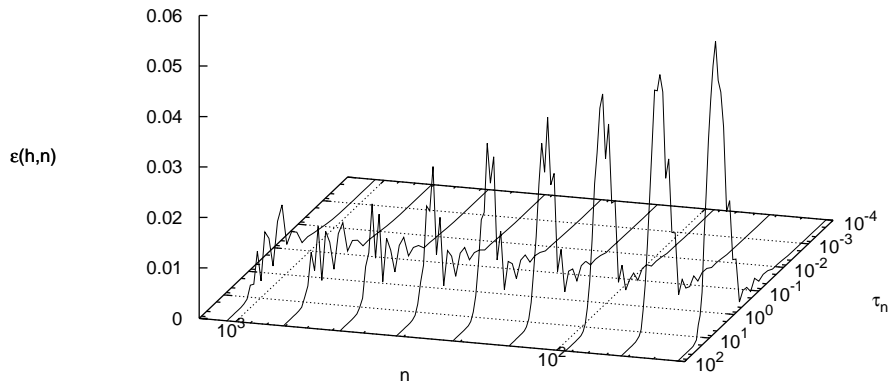}}
\caption{From Figure \ref{figp25b}: $\varepsilon(h,n)=|A(h,n)- e^{-\theta \tau_n}|$ for the M\"obius dynamics $\varphi_3$ and invariant Minkowski measure $d\mu=dQ$ versus $n$ and $\tau_n$. }
\label{figp25c}
\end{figure}

We can now perform the usual analysis, employing the the scaling induced by eq. (\ref{eq-mink41}). Results are shown in Figure \ref{figp25b}. Convergence to the asymptotic extremal law is again verified, but we observe numerically the presence of a non--trivial extremal exponent $\theta \simeq 0.47$. We leave the rigorous derivation of this exponent to future work, but we remark that its origin must be traced back to the fact that, since $\varphi_3(1/k) = 1/(k-1)$, the set $K$ is a trajectory of the motion.  Figure \ref{figp25c} plots as before differences of experimental data and asymptotic theory.

\section{Conclusions}
In this paper we have introduced a new paradigm in the dynamical theory of extreme events: the r\^ole of fractal intensity functions, lacking any convexity or maximality property.

We started by describing a toy model, the Cantor Ladder, that was then extended to general intensity functions with singularities on fractal subsets of the dynamical attractor. We have shown that the conventional theory can be generalized to these cases. The related extreme value laws bring into play familiar fractal quantities, the Minkowski dimension and content, as well as they generalizations, defined herein, to singular continuous dynamical measures. Remarkably, these laws appear to hold also for ergodic, non--mixing systems.
We have finally introduced the notion of extremely rare event and the corresponding extreme value laws, exemplified by a dynamical system with origin in number theory.

We hope that the theoretical framework and numerical results presented in this paper will stimulate further analysis---to prove rigorously the facts demonstrated by compelling numerical evidence, to clarify how complexity of the fractal landscape of the intensity function plays a r\^ole in guaranteeing the validity of extreme value laws---and at the same time be useful for applications of the theory to real phenomena.

\section{Acknowledgements}
G. M. acknowledges hospitality at the Isaac Newton Institute for Mathematical Sciences (Cambridge, UK) during the 2013 programme {\em Mathematics for the Fluid Earth}, where this work has been started. Generous support has also been provided by ANR Perturbations, ANR Valet and AMIDEX for participation of GM to the workshop {\em Dynamiques \`a Porquerolles 2015}. Fruitful conversations with S. Vaienti and J. Freitas and the other participants to the two programs, whose names are too many to list, are also kindly acknowledged.
Computations for this paper were performed at the INFN computer center at Pisa.

\end{document}